\documentclass[10pt,reqno]{amsart}
\usepackage{amsthm}
\usepackage{amsmath,amssymb}
\usepackage{harpoon}
\usepackage{threeparttable}
\usepackage{tabularx,multirow,bigdelim}
\usepackage{algorithm,algorithmic}
\usepackage{graphicx,subfigure}
\usepackage{hyperref}

\usepackage{mathrsfs}

\theoremstyle{plain}
\newtheorem{theorem}{Theorem}[section]
\newtheorem{lemma}{Lemma}[section]

\theoremstyle{definition}
\newtheorem{definition}{Definition}[section]

\theoremstyle{remark}

\newenvironment{prof}[1][Proof]{\indent\textbf{#1}\quad}{\hfill $\Box$\vspace{0.7mm}}
\numberwithin{equation}{section}

\def\no{\nonumber}
\def\R{\mathbb{R}}

\allowdisplaybreaks[4]

\begin{document}

\title{A note on the first variation of the total mass}

\author{Xiaomei Sun}
\address{College of Informatics, Huazhong Agricultural University, Wuhan 430070, China}
\email{xmsunn@mail.hzau.edu.cn}

\author{Anqiang Zhu}
\address{School of Mathematics and Statistics, Wuhan University, Wuhan 430072, China}
\email{aqzhu.math@whu.edu.cn}
\thanks{The second author is supported by the Hubei Provincial Natural Science Foundation of China (Grant No. 2025AFB811).}

\begin{abstract}
In this paper, we establish a proof for the first variation formula of the total mass within the $L_p$ framework. Our main result removes an extra restrictive determinant condition imposed in a theorem originally proved in \cite{Yedeping}.
\end{abstract}


\maketitle

\section{Introduction}
This work is concerned with the first variation of the total mass functional associated with log-concave functions. For a log-concave function $f=e^{-\varphi}$, where $\varphi:\R^n\to\R\cup\{+\infty\}$ is convex and satisfies $\lim_{|x|\to\infty}\varphi(x)=+\infty$, we define its total mass as follows.

\begin{definition}
The total mass functional $J$ acting on log-concave functions is given by
\begin{align}
J(f):=\int_{\R^n}f(x)\,dx.
\end{align}
\end{definition}

We focus on variations with respect to the $L_p$ Asplund sum (see Definition \ref{def:pAsplund}). The first variation of $J$ at $f$ along a $p$-perturbation $g$ is defined below.

\begin{definition}
For $p>1$ and $f,g\in\mathscr{A}_0$, the $L_p$-first variation of $J$ is
\begin{align}
\delta J_p(f,g):=\lim_{t\to0^+}\frac{J(f\oplus_p t\cdot_p g)-J(f)}{t}.
\end{align}
The function class $\mathscr{A}_0$ will be specified in Section \ref{sec:nota}.
\end{definition}

The variational theory for log-concave functions was initiated by Colesanti and Fragalà, who derived the first variation of total mass for the classical $p=1$ Asplund sum in their foundational work \cite{Colesanti}. Later, Rotem refined their variational formula and argued that essential continuity constitutes the minimal, optimal set of assumptions guaranteeing validity of such variation identities \cite{Rotem}.

For the range $p>1$, Fang, Niufa, Xing and Ye extended the $p=1$ results to the $L_p$ setting in \cite{Yedeping}, proving the following variational formula under an additional determinant compatibility condition.

\begin{theorem}[\cite{Yedeping}]\label{thm:oldvar}
Let $f=e^{-\varphi}\in\mathscr{A}_0'$ and $g=e^{-\psi}\in\mathscr{A}_0'$ with $p>1$, and suppose $g$ is an admissible $p$-perturbation of $f$. Assume further that there exists a constant $k>0$ such that
\begin{equation}\label{eq:compcond}
\det\big(\nabla^2(\varphi^*)^p(y)\big)\leq k \cdot \big(\varphi^*(y)\big)^{n(p-1)} \det\big(\nabla^2\varphi^*(y)\big)
\end{equation}
holds for all $y\in\R^n\setminus\{0\}$. Then
\begin{align}
\delta J_p(f,g)
&=\frac{1}{p}\int_{\R^n}\big(\psi^*(\nabla\varphi(x))\big)^p \big(\varphi^*(\nabla\varphi(x))\big)^{1-p}e^{-\varphi(x)}\,dx \no\\
&=\frac{1}{p}\int_{\R^n}\big(\psi^*(y)\big)^p \big(\varphi^*(y)\big)^{1-p}\,d\mu(f,y).\label{eq:oldformula}
\end{align}
\end{theorem}

The classes $\mathscr{A}_0'$ and admissible $p$-perturbations are defined in Section \ref{sec:nota}. When $p=1$, formula \eqref{eq:oldformula} recovers the identity from \cite{Colesanti}, yet the extra compatibility constraint \eqref{eq:compcond} is entirely absent from the $p=1$ theory. 

In this short note, we eliminate the restrictive determinant condition \eqref{eq:compcond} for all $p>1$ and prove the following simplified main theorem.

\begin{theorem}\label{thm:main}
Let $f=e^{-\varphi}\in\mathscr{A}_0'$, $g=e^{-\psi}\in\mathscr{A}_0'$ with $p>1$, and let $g$ be an admissible $p$-perturbation of $f$. Then
\begin{align}
\delta J_p(f,g)
&=\frac{1}{p}\int_{\R^n}\big(\psi^*(\nabla\varphi(x))\big)^p \big(\varphi^*(\nabla\varphi(x))\big)^{1-p}e^{-\varphi(x)}\,dx \no\\
&=\frac{1}{p}\int_{\R^n}\big(\psi^*(y)\big)^p \big(\varphi^*(y)\big)^{1-p}\,d\mu(f,y).\label{eq:mainformula}
\end{align}
\end{theorem}

\section{Preliminaries and Notation}\label{sec:nota}
We first introduce standard classes of convex functions and associated log-concave functions used throughout the paper. Let $\mathcal{C}$ denote the collection of all proper convex functions $\R^n\to\R\cup\{+\infty\}$. Define
\[
\mathscr{L}=\big\{\varphi\in\mathcal{C}\,\big|\, \lim_{|x|\to+\infty}\varphi(x)=+\infty\big\},\quad
\mathscr{A}=\{f=e^{-\varphi}\mid \varphi\in\mathscr{L}\}.
\]
We restrict attention to nonnegative convex functions vanishing at the origin:
\[
\mathscr{L}_0=\big\{\varphi\in\mathscr{L}\,\big|\, \varphi\geq0,\ \varphi(0)=0,\ \varphi\text{ lower semicontinuous}\big\},\quad
\mathscr{A}_0=\{e^{-\varphi}\mid \varphi\in\mathscr{L}_0\}.
\]

Our analysis focuses on a smooth subclass $\mathscr{L}_0'\subset\mathscr{L}_0$ consisting of strictly convex supercoercive functions with full domain:
\[
\begin{split}
\mathscr{L}_0':=\big\{
\varphi\in\mathscr{L}_0\,\big|\,
&\mathrm{dom}(\varphi)=\R^n,\ \varphi\in C^1(\R^n)\cap C_+^2(\R^n\setminus\{0\}),\\
&\varphi\text{ strictly convex, supercoercive}
\big\}.
\end{split}
\]
Here $C_{+}^{2}(\mathbb{R}^{n}\setminus\{0\})$ denotes the set of functions such that every function $\varphi$ restricted to $\mathbb{R}^{n}\setminus\{0\}$ is twice continuously differentiable, and its Hessian matrix $\nabla^{2}\varphi(x)$ is positive definite at each point $x\in\mathbb{R}^{n}\setminus\{0\}$.

A convex function $\varphi$ is called supercoercive if $\lim_{|x|\to\infty}\varphi(x)/|x|=+\infty$. The corresponding log-concave class is $\mathscr{A}_0'=\{e^{-\varphi}\mid \varphi\in\mathscr{L}_0'\}$.

We now recall the definition of the $L_p$ Asplund sum for log-concave functions, following \cite{Yedeping}.

\begin{definition}[\cite{Yedeping}]\label{def:pAsplund}
Take $p>1$, $f=e^{-\varphi}\in\mathscr{A}_0$, $g=e^{-\psi}\in\mathscr{A}_0$. The $L_p$ Asplund sum $f\oplus_p g$ is defined via the convex conjugate operation as $f\oplus_p g=e^{-\varphi\Box_p\psi}$, where
\[
\varphi\Box_p\psi=\left[\big((\varphi^*)^p + (\psi^*)^p\big)^{1/p}\right]^*.
\]
Here $\varphi^*(y)=\sup_{x\in\R^n}\big(\langle x,y\rangle-\varphi(x)\big)$ denotes the Legendre-Fenchel conjugate of $\varphi$.

For $\alpha>0$, the $p$-scaling of a convex function is
\[
(\varphi\cdot_p \alpha)(x)=\alpha^{1/p}\varphi\big(\alpha^{-1/p}x\big),\quad \forall x\in\R^n.
\]
For $\alpha,\beta>0$, the scaled $L_p$ Asplund sum satisfies
\[
\alpha\cdot_p f\oplus_p \beta\cdot_p g = e^{-\varphi\cdot_p\alpha \Box_p \psi\cdot_p\beta},\quad
\varphi\cdot_p\alpha \Box_p \psi\cdot_p\beta
=\left[\big(\alpha(\varphi^*)^p+\beta(\psi^*)^p\big)^{1/p}\right]^*.
\]
\end{definition}

The following monotonicity property of the $L_p$ Asplund sum will be critical for uniform convergence arguments.

\begin{lemma}[\cite{Yedeping}]\label{lem:monotone}
Let $p>1$, $t>0$, and set $\varphi_t=\varphi\Box_p(\psi\cdot_p t)$, $f_t=e^{-\varphi_t}$. For any $0<s<t\leq1$ and all $x\in\R^n$,
\[
\varphi_1(x)\leq\varphi_t(x)\leq\varphi_s(x)\leq\varphi(x),\quad
f(x)\leq f_s(x)\leq f_t(x)\leq f_1(x).
\]
\end{lemma}

\begin{definition}
A pair $(C,f)$ with open convex domain $C=\mathrm{int}\,\mathrm{dom}\,f$ is said to be of Legendre type if:
\begin{enumerate}
    \item $C\neq\emptyset$;
    \item $f$ is differentiable everywhere on $C$;
    \item For every sequence $\{x_i\}\subset C$ converging to a boundary point of $C$, one has $\lim_{i\to\infty}|\nabla f(x_i)|=+\infty$.
\end{enumerate}
\end{definition}

Given a differentiable function $f$ on an open convex set $C$, its Legendre conjugate pair $(D,g)$ is defined as follows: let $D=\nabla f(C)$, and set
\[
g(y)=\langle x,y\rangle-f(x),\quad y\in D,\ x=(\nabla f)^{-1}(y).
\]

\begin{lemma}[\cite{Rockafellar,Yedeping}]\label{lem:legendredual}
Suppose $\phi\in C^1(\R^n)$ and $\mathrm{dom}(\phi^*)=\R^n$. Then $(\R^n,\phi)$ is Legendre type if and only if $(\R^n,\phi^*)$ is Legendre type. In this case, $\nabla\phi:\R^n\to\R^n$ and $\nabla\phi^*:\R^n\to\R^n$ are mutually inverse continuous bijections, i.e., $(\nabla\phi)^{-1}=\nabla\phi^*$.
\end{lemma}

\begin{lemma}[\cite{Yedeping}]\label{lem:regconj}
If $\varphi\in\mathscr{L}_0'$, then its conjugate satisfies $\varphi^*\in\mathscr{L}_0'$. Moreover, $\nabla\varphi(0)=\nabla\varphi^*(0)=0$, and
\[
\{x\in\R^n\mid \varphi(x)=0\}=\{x\in\R^n\mid \varphi^*(x)=0\}=\{0\}.
\]
\end{lemma}

\begin{lemma}[\cite{Yedeping}]\label{lem:phitleg}
Let $\varphi,\psi\in\mathscr{L}_0'$, $p>1$, $t>0$, and define $\varphi_t=\varphi\Box_p(\psi\cdot_p t)$. Then $\varphi_t\in\mathscr{L}_0'$, and both $(\R^n,\varphi_t)$ and $(\R^n,\varphi_t^*)$ are Legendre-type convex functions.
\end{lemma}

\begin{definition}[\cite{Yedeping}]
Let $p>1$, $f=e^{-\varphi}\in\mathscr{A}_0$. We call $g=e^{-\psi}\in\mathscr{A}_0$ an admissible $p$-perturbation of $f$ if there exists a constant $c>0$ such that $(\varphi^*)^p-c(\psi^*)^p$ is convex on $\R^n$.
\end{definition}

Convexity of $(\varphi^*)^p-c(\psi^*)^p$ yields a global lower bound at the origin:
\begin{align}\label{eq:convexlowbd}
\big((\varphi^*)^p-c(\psi^*)^p\big)(y)
\geq \big((\varphi^*)^p-c(\psi^*)^p\big)(0)
+\big\langle y,\nabla\big((\varphi^*)^p-c(\psi^*)^p\big)(0)\big\rangle=0,\quad \forall y\in\R^n.
\end{align}

The time derivative of $\varphi_t$ is computed explicitly in the next lemma.

\begin{lemma}[\cite{Yedeping}]\label{lem:dphitdt}
For $\varphi,\psi\in\mathscr{L}_0'$, $p>1$, $t>0$, set $\varphi_t=\varphi\Box_p(\psi\cdot_p t)$. Then
\begin{align}
\frac{d}{dt}\varphi_t(x)=
\begin{cases}
-\displaystyle\frac{1}{p}\big(\psi^*(\nabla\varphi_t(x))\big)^p \big(\varphi_t^*(\nabla\varphi_t(x))\big)^{1-p}, & x\neq0,\\[4pt]
0, & x=0.
\end{cases}
\end{align}
In particular, $\partial_t\varphi_t(0)\big|_{t=0^+}=0$, and for all $x\neq0$,
\begin{align}\label{eq:derivatt0}
\left.\frac{d}{dt}\varphi_t(x)\right|_{t=0^+}
=-\frac{1}{p}\big(\psi^*(\nabla\varphi(x))\big)^p \big(\varphi^*(\nabla\varphi(x))\big)^{1-p}.
\end{align}
\end{lemma}

All functions in $\mathscr{L}$ enjoy linear lower growth:
\begin{lemma}[\cite{Colesanti}]\label{lem:linearlower}
For any $u\in\mathscr{L}$, there exist constants $a>0$, $b\in\R$ such that
\[
u(x)\geq a|x|+b,\quad \forall x\in\R^n.
\]
Its conjugate $u^*$ is proper and bounded below everywhere on $\R^n$.
\end{lemma}

We will also rely on uniform convergence of gradients for pointwise convergent convex functions:
\begin{theorem}[\cite{Rockafellar}]\label{thm:gradunifconv}
Let $C\subset\R^n$ be open convex, and $f,f_i$ convex, finite and differentiable on $C$. If $f_i(x)\to f(x)$ pointwise for all $x\in C$, then $\nabla f_i\to\nabla f$ uniformly on every compact subset of $C$.
\end{theorem}

\section{Proof of Theorem \ref{thm:main}}
We establish several auxiliary continuity and uniform convergence lemmas before proving the main variational identity.

\begin{lemma}\label{lem:phitunifconv}
Fix $s\geq0$. For every compact set $E\subset\R^n$, $\varphi_t(x)\to\varphi_s(x)$ uniformly in $x\in E$ as $t\to s$.
\end{lemma}
\begin{prof}
By monotonicity (Lemma \ref{lem:monotone}), $\varphi_t(x)$ decreases in $t\in(0,1]$. For $0<s<t\leq1$, $\varphi_t(x)\leq\varphi_s(x)$.
Fix compact $E\subset\R^n$. Since $\varphi_s\in\mathscr{L}_0'$, $\nabla\varphi_s$ is continuous, so $\nabla\varphi_s(E)$ is contained in some closed ball $B_r$. The map $z\mapsto z^{1/p}$ is concave on $[0,\infty)$ with $0^{1/p}=0$, hence subadditive: $(A+B)^{1/p}\leq A^{1/p}+B^{1/p}$ for $A,B\geq0$. Write
\[
(\varphi_t^*)^p = (\varphi^*)^p + s(\psi^*)^p + (t-s)(\psi^*)^p,
\]
so
\[
\varphi_t^*=\big((\varphi^*)^p+s(\psi^*)^p+(t-s)(\psi^*)^p\big)^{1/p}
\leq \varphi_s^* + (t-s)^{1/p}\psi^*.
\]
For any $x\in E$,
\begin{align*}
\varphi_s(x)
&\geq \varphi_t(x)=\sup_{y\in\R^n}\big(\langle x,y\rangle-\varphi_t^*(y)\big)
\geq \sup_{y\in B_r}\big(\langle x,y\rangle-\varphi_t^*(y)\big)\\
&\geq \sup_{y\in B_r}\big(\langle x,y\rangle-\varphi_s^*(y)\big)-(t-s)^{1/p}\sup_{y\in B_r}\psi^*(y)\\
&=\big\langle x,\nabla\varphi_s(x)\big\rangle-\varphi_s^*(\nabla\varphi_s(x)) - C(t-s)^{1/p}
=\varphi_s(x)-C(t-s)^{1/p},
\end{align*}
where $C=\sup_{B_r}\psi^*<\infty$. Thus
\[
|\varphi_t(x)-\varphi_s(x)|\leq C\cdot |t-s|^{1/p},\quad \forall x\in E,
\]
which yields uniform convergence on compact sets. The case $s=0$ follows identically, giving uniform convergence $\varphi_t\to\varphi$ as $t\to0^+$.
\end{prof}

Combined with Theorem \ref{thm:gradunifconv}, $\nabla\varphi_t(x)\to\nabla\varphi_s(x)$ uniformly on all closed balls as $t\to s$.

\begin{lemma}\label{lem:phitjointcont}
The map $(x,t)\mapsto\varphi_t(x)$ is continuous on $\R^n\times [0,\infty)$.
\end{lemma}
\begin{prof}
Take arbitrary $(s_0,x_0)\in[0,\infty)\times\R^n$. Since $\varphi_{s_0}\in\mathscr{L}_0'$, it is continuous in $x$, so there exists $\delta_0>0$ such that $|\varphi_{s_0}(x)-\varphi_{s_0}(x_0)|<\varepsilon$ whenever $|x-x_0|<\delta_0$.
By Lemma \ref{lem:phitunifconv}, $\varphi_t$ converges uniformly to $\varphi_{s_0}$ on the closed ball $\overline{B(x_0,\delta_0)}$ as $t\to s_0$. Hence there exists $\delta_1>0$ such that $|\varphi_t(x)-\varphi_{s_0}(x)|<\varepsilon$ for all $|t-s_0|<\delta_1$ and $x\in\overline{B(x_0,\delta_0)}$. Set $\delta=\min\{\delta_0,\delta_1\}$. If $|x-x_0|+|t-s_0|<\delta$, then
\[
|\varphi_t(x)-\varphi_{s_0}(x_0)|
\leq |\varphi_t(x)-\varphi_{s_0}(x)| + |\varphi_{s_0}(x)-\varphi_{s_0}(x_0)| < 2\varepsilon.
\]
Joint continuity follows.
\end{prof}

As an immediate corollary, $f_t(x)=e^{-\varphi_t(x)}$ is continuous on $[0,\infty)\times\R^n$.

Define the integrand 
\[
F(s,x):=\big(\psi^*(\nabla\varphi_s(x))\big)^p \big(\varphi_s^*(\nabla\varphi_s(x))\big)^{1-p}.
\]

\begin{lemma}\label{lem:Fcontaway0}
$F(s,x)$ is continuous on $[0,\infty)\times(\R^n\setminus\{0\})$.
\end{lemma}
\begin{prof}
We first prove joint continuity of $(s,x)\mapsto\nabla\varphi_s(x)$. Fix $(s_0,x_0)$. Since $\nabla\varphi_{s_0}$ is continuous in $x$, there exists $\delta_1>0$ such that $|\nabla\varphi_{s_0}(x)-\nabla\varphi_{s_0}(x_0)|<\varepsilon$ for $|x-x_0|<\delta_1$. By Theorem \ref{thm:gradunifconv}, $\nabla\varphi_s$ converges uniformly to $\nabla\varphi_{s_0}$ on $\overline{B(x_0,\delta_1)}$, so there exists $\delta_2>0$ such that $|\nabla\varphi_s(x)-\nabla\varphi_{s_0}(x)|<\varepsilon$ whenever $|s-s_0|<\delta_2$ and $x\in\overline{B(x_0,\delta_1)}$. For $|x-x_0|+|s-s_0|<\min\{\delta_1,\delta_2\}$,
\[
|\nabla\varphi_s(x)-\nabla\varphi_{s_0}(x_0)|
\leq |\nabla\varphi_s(x)-\nabla\varphi_{s_0}(x)| + |\nabla\varphi_{s_0}(x)-\nabla\varphi_{s_0}(x_0)| < 2\varepsilon.
\]
Thus $\nabla\varphi_s(x)$ is jointly continuous.
The functions $\varphi^*,\psi^*$ belong to $\mathscr{L}_0'$ (Lemma \ref{lem:regconj}), hence continuous, and $\nabla\varphi_s(x)\neq0$ for all $x\neq0$. The expression defining $F(s,x)$ is a composition of continuous maps away from $\{0\}$, which yields the claimed continuity.
\end{prof}

\begin{lemma}\label{lem:Fvanishatorigin}
For any fixed $s_0\in[0,\infty)$, $\lim_{(t,x)\to(s_0,0)}F(t,x)=0$.
\end{lemma}
\begin{prof}
Admissible $p$-perturbation condition supplies $c>0$ with $(\varphi^*)^p(y)\geq c(\psi^*)^p(y)$ for all $y$, so
\[
\frac{(\psi^*)^p(\nabla\varphi_t(x))}{(\varphi^*)^p(\nabla\varphi_t(x))+t(\psi^*)^p(\nabla\varphi_t(x))}\leq \frac{1}{c}.
\]
We bound
\[
0\leq F(t,x)
\leq \frac{1}{c}\big((\varphi^*)^p(\nabla\varphi_t(x))+t(\psi^*)^p(\nabla\varphi_t(x))\big)^{1/p}.
\]
Pass to the limit $(t,x)\to(s_0,0)$ and use continuity of $\nabla\varphi_t,\varphi^*,\psi^*$:
\[
\lim_{(t,x)\to(s_0,0)}\big((\varphi^*)^p(\nabla\varphi_t(x))+t(\psi^*)^p(\nabla\varphi_t(x))\big)^{1/p}
=\big((\varphi^*)^p(0)+s_0(\psi^*)^p(0)\big)^{1/p}=0,
\]
since $\varphi^*(0)=\psi^*(0)=0$. The squeeze argument gives $\lim_{(t,x)\to(s_0,0)}F(t,x)=0$.
\end{prof}

Combining Lemmas \ref{lem:Fcontaway0} and \ref{lem:Fvanishatorigin} with Lemma \ref{lem:dphitdt}, the partial derivative $\partial_t\varphi_t(x)$ is continuous on $[0,1]\times\R^n$.

\begin{lemma}\label{lem:intunifconv}
The integral
\[
I(s)=\frac{1}{p}\int_{\R^n}e^{-\varphi_s(x)} F(s,x)\,dx
\]
converges uniformly for $s\in[0,1]$.
\end{lemma}
\begin{prof}
From \eqref{eq:convexlowbd}, $(\psi^*)^p\leq c^{-1}(\varphi^*)^p$, so
\[
F(s,x)\leq \frac{1}{c}\big(\varphi_s^*(\nabla\varphi_s(x))\big).
\]
It suffices to verify uniform convergence of
\[
\tilde{I}(s)=\frac{1}{pc}\int_{\R^n}e^{-\varphi_s(x)}\varphi_s^*(\nabla\varphi_s(x))\,dx.
\]
Recall the Legendre identity $\varphi_s^*(\nabla\varphi_s(x))=\langle x,\nabla\varphi_s(x)\rangle-\varphi_s(x)$. Integrate over $\R^n\setminus B_R$:
\begin{align*}
\int_{\R^n\setminus B_R}e^{-\varphi_s(x)}\varphi_s^*(\nabla\varphi_s(x))dx
&=\int_{\R^n\setminus B_R}\big(\langle x,\nabla\varphi_s(x)\rangle-\varphi_s(x)\big)e^{-\varphi_s(x)}dx\\
&=-\int_{\R^n\setminus B_R}\langle x,\nabla e^{-\varphi_s(x)}\rangle dx
-\int_{\R^n\setminus B_R}\varphi_s(x)e^{-\varphi_s(x)}dx.
\end{align*}
Apply the divergence theorem to the vector field $xe^{-\varphi_s(x)}$:
\[
\int_{\R^n\setminus B_R}\langle x,\nabla e^{-\varphi_s(x)}\rangle dx
= -\int_{\partial B_R} R e^{-\varphi_s(x)}d\mathcal{H}^{n-1}
-n\int_{\R^n\setminus B_R}e^{-\varphi_s(x)}dx.
\]
Thus
\begin{align*}
\int_{\R^n\setminus B_R}e^{-\varphi_s(x)}\varphi_s^*(\nabla\varphi_s(x))dx
&=\int_{\partial B_R}R e^{-\varphi_s}d\mathcal{H}^{n-1}
+n\int_{\R^n\setminus B_R}e^{-\varphi_s}dx
-\int_{\R^n\setminus B_R}\varphi_s e^{-\varphi_s}dx.
\end{align*}
From admissibility, $(\varphi_s^*)^p=(\varphi^*)^p+s(\psi^*)^p\leq (1+s/c)(\varphi^*)^p$, so
\[
\varphi_s(x)\geq \left(\frac{c+s}{c}\right)^{1/p}\varphi\left(\left(\frac{c+s}{c}\right)^{-1/p}x\right).
\]
Using the linear lower bound Lemma \ref{lem:linearlower}:
\[
\varphi_s(x)\geq a|x|+b,\quad \forall s\in[0,1],\ x\in\R^n.
\]
This uniform linear lower bound yields three uniform decay estimates as $R\to\infty$:
\[
\int_{\partial B_R}R e^{-\varphi_s}d\mathcal{H}^{n-1}
\leq R e^{-aR-b}\mathcal{H}^{n-1}(\partial B_R)\to0,\quad
\int_{\R^n\setminus B_R}e^{-\varphi_s}dx\leq\int_{|x|>R}e^{-a|x|-b}dx\to0,
\]
and for sufficiently large $R$ with $a|x|+b>1$ on $|x|>R$, the map $z\mapsto z e^{-z}$ is decreasing for $z>1$, hence
\[
\int_{\R^n\setminus B_R}\varphi_s e^{-\varphi_s}dx
\leq \int_{|x|>R}(a|x|+b)e^{-a|x|-b}dx\to0,
\]
all convergences uniform in $s\in[0,1]$. This completes uniform integrability.
\end{prof}

\begin{prof}[Proof of Theorem \ref{thm:main}]
By definition of the $L_p$-first variation:
\[
\delta J_p(f,g)=\lim_{t\to0^+}\frac{J(f\oplus_p t\cdot_p g)-J(f)}{t}
=\lim_{t\to0^+}\int_{\R^n}\frac{e^{-\varphi_t(x)}-e^{-\varphi(x)}}{t}dx.
\]
The function $f(x,t)=e^{-\varphi_t(x)}$ is jointly continuous on $\R^n\times [0,\infty)$ (Lemma \ref{lem:phitjointcont}), and each integral $J(f_t)=\int_{\R^n}e^{-\varphi_t}dx$ is finite via Lemma \ref{lem:linearlower}.
The partial derivative 
\begin{align}
    \partial_t e^{-\varphi_t(x)}=-e^{-\varphi_t(x)}\partial_t\varphi_t(x)=e^{\varphi_{t}(x)}\frac{1}{p}F(t,x)
\end{align}
 is continuous on $\R^n\times[0,1]$ (Lemma
  \ref{lem:Fcontaway0}, Lemma \ref{lem:Fvanishatorigin}) , and $\int_{\R^n}\partial_t e^{-\varphi_t}dx$ converges uniformly for $t\in[0,1]$ (Lemma \ref{lem:intunifconv}). We apply the parameter differentiation theorem (Theorem \ref{thm:diffunderint}) from the appendix to interchange limit and integral:
\begin{align*}
\delta J_p(f,g)
&=\lim_{t\to0^+}\int_{\R^n}\frac{e^{-\varphi_t(x)}-e^{-\varphi(x)}}{t}dx
=\int_{\R^n}\left.\frac{\partial}{\partial t}e^{-\varphi_t(x)}\right|_{t=0^+}dx\\
&=\frac{1}{p}\int_{\R^n}\big(\psi^*(\nabla\varphi(x))\big)^p \big(\varphi^*(\nabla\varphi(x))\big)^{1-p}e^{-\varphi(x)}dx.
\end{align*}
The dual measure identity follows from the standard change of variables via the Legendre gradient bijection $\nabla\varphi$.
\end{prof}

\section{Appendix: Differentiation Under the Integral Sign}
We include two standard integral differentiation theorems for self-containedness of the paper.

\begin{theorem}\label{thm:compactdiff}
Let $K\subset\R^n$ compact, $f:K\times[a,b]\to\R$ continuous in $x$ for each fixed $y$, with continuous partial derivative $\partial_y f\in C(K\times[a,b])$. Then
\[
\frac{d}{dy}\int_K f(x,y)dx=\int_K \partial_y f(x,y)dx,\quad \forall y\in[a,b].
\]
\end{theorem}
\begin{prof}
Fix $y_0\in[a,b]$, set $I(y)=\int_K f(x,y)dx$. For small $|k|$, the mean value theorem gives
\[
\frac{f(x,y_0+k)-f(x,y_0)}{k}=\partial_y f(x,y_0+\theta k),\quad \theta\in(0,1).
\]
Since $\partial_y f$ is uniformly continuous on the compact product $K\times[a,b]$, for any $\varepsilon>0$ there exists $\delta>0$ such that $|\partial_y f(x,y_0+\theta k)-\partial_y f(x,y_0)|<\varepsilon$ whenever $|k|<\delta$. Thus
\[
\left|\int_K\frac{f(x,y_0+k)-f(x,y_0)}{k}dx-\int_K\partial_y f(x,y_0)dx\right|
\leq \varepsilon\cdot\mathrm{Vol}(K),
\]
and the limit as $k\to0$ yields the desired equality.
\end{prof}

\begin{theorem}[\cite{Rudin-baby}]\label{thm:diffunderint}
Suppose $f(x,y):\R^n\times[a,b]\to\R, x\in \R^{n}, y\in [a,b]$ is continuous, $\partial_y f$ continuous everywhere, and $I(y)=\int_{\R^n}f(x,y)dx$ exists for all $y\in[a,b]$. If $\int_{\R^n}\partial_y f(x,y)dx$ converges uniformly over $y\in[a,b]$, then
\[
I'(y)=\int_{\R^n}\partial_y f(x,y)dx,\quad \forall y\in[a,b].
\]
\end{theorem}
\begin{prof}
Fix $y_0\in[a,b]$, write $B_R=B(0,R)$. Define
\[
F(R,k)=\int_{B_R}\frac{f(x,y_0+k)-f(x,y_0)}{k}dx,\quad
\varphi(k)=\int_{\R^n}\frac{f(x,y_0+k)-f(x,y_0)}{k}dx.
\]
By Theorem \ref{thm:compactdiff}, $\lim_{k\to0}F(R,k)=\int_{B_R}\partial_y f(x,y_0)dx$ for every fixed $R>0$.
Uniform convergence of $\int_{\R^{n}}\partial_y f(x,y) dx$ implies that for all $\varepsilon>0$, there exists $A_0>0$ such that
\[
\left|\int_{\overline{B_{R_1}}\setminus B_{R_2}}\partial_y f(x,y)dx\right|<\varepsilon,\quad \forall y\in[a,b],\ R_1>R_2>A_0.
\]
Here $B_{R}$ are open ball.
Let $\eta(y)=\int_{\overline{B_{R_1}}\setminus B_{R_2}}f(x,y)dx$. By Theorem \ref{thm:compactdiff}, we have 
\begin{align*}
    |\eta'(y)|=\left|\int_{\overline{B_{R_1}}\setminus B_{R_2}}\partial_y f(x,y)dx\right|\leq \varepsilon.
\end{align*}
 for all $y\in[a,b]$. Then the mean value theorem gives
\[
\left|\int_{\overline{B_{R_1}}\setminus B_{R_2}}\frac{f(x,y_0+k)-f(x,y_0)}{k}dx\right|=|\eta'(y_0+\theta k)|<\varepsilon.
\]
The above inequality is valid for any $R_{1}$ sufficiently large.
Choose $R>A_0$ large enough so that both
\begin{align*}
    |\varphi(k)-F(R,k)|=\left|\int_{\R^{n}\setminus B_{R}} \frac{f(x,y_0+k)-f(x,y_0)}{k}\right| \leq \varepsilon
\end{align*}
and 
\begin{align*}
  \left|\int_{\R^n\setminus B_R}\partial_y f(x,y_0)dx\right|<\varepsilon.  
\end{align*}
For this fixed $R$, pick $\delta>0$ with $|F(R,k)-\int_{B_R}\partial_y f(x,y_0)dx|<\varepsilon$ whenever $|k|<\delta$. The triangle inequality yields
\begin{align*}
    \left|\varphi(k)-\int_{\R^n}\partial_y f(x,y_0)dx\right|
&\leq |\varphi(k)-F(R,k)| + |F(R,k)-\int_{B_R}\partial_y f(x,y_{0})|dx \\
&+ \left|\int_{\R^n\setminus B_R}\partial_y f(x,y_{0})\right|dx
<3\varepsilon.
\end{align*}

Let $\varepsilon\to0$ to conclude $\lim_{k\to0}\varphi(k)=\int_{\R^n}\partial_y f(x,y_0)dx$.
\end{prof}

\section*{Declarations}
\subsection*{Author Contributions}
Xiaomei Sun and Anqiang Zhu contributed equally to this manuscript. Both authors derived the results, wrote the draft, revised all proofs, and approved the final version.

\subsection*{Conflict of Interest}
The authors declare no competing interests.

\subsection*{Data Availability}
This is a purely theoretical work; no experimental data was generated or analysed.

\bibliographystyle{amsplain}
\bibliography{Ref.bib}
\end{document}